\documentclass[11pt]{article}

\usepackage{latexsym}
\usepackage{amssymb}

\newtheorem{theorem}{Theorem}

\newtheorem{definition}{Definition}

\begin{document}
{
\begin{center}
{\Large\bf
On the generalized hypergeometric function, Sobolev orthogonal polynomials and biorthogonal rational functions.}
\end{center}
\begin{center}
{\bf S.M. Zagorodnyuk}
\end{center}

\section{Introduction.}

The theories of orthogonal polynomials on the real line and on the unit circle have a lot of contributions
and applications~\cite{cit_50000_Gabor_Szego}, \cite{cit_5000_Ismail}, \cite{cit_48000_Simon_1}, \cite{cit_48000_Simon_2}.
One of their possible generalizations, the theory of Sobolev orthogonal polynomials is now studied 
intensively by many mathematicians~(see a survey in~\cite{cit_5150_M_X}). 
Important ingredients, which supported the importance of classical systems of polynomials $\{ p_n(z) \}_{n=0}^\infty$, are the recurrence
relation and the differential equations for $p_n(z)$. Thus, it is natural to seek for such properties of Sobolev
orthogonal polynomials.
In this paper we shall provide a large class of polynomials $g_n(z)$ which has both these properties. Moreover,
polynomials $g_n(z)$ are related to biorthogonal rational functions and Jacobi-type pencils.

Let $p,q$ be some fixed non-negative integers. Denote
$$ g_n(z) = g_n(z;a_1,...,a_p; b_1,...,b_q) = $$
\begin{equation}
\label{f1_5} 
= \sum_{k=0}^n \frac{(a_1)_k ... (a_p)_k}{(b_1)_k ... (b_q)_k} \frac{z^k}{k!},\qquad n=0,1,2,..., 
\end{equation}
where $a_j,b_l\in\mathbb{C}\backslash\{ 0,-1,-2,... \}$.
Thus, $g_n(z)$ is the $n$-th partial sum of
the generalized hypergeometric series ${}_p F_q(a_1,...,a_p; b_1,...,b_q;z)$, and $\deg g_n=n$.
As usual, $p=0$ ($q=0$) means that $(a_j)_k$ (respectively $(b_l)_k)$ are absent.
By $G_n(z)$ we denote the corresponding monic polynomials:
\begin{equation}
\label{f1_7} 
G_n(z) = \frac{ n! (b_1)_n ... (b_q)_n }{ (a_1)_n ... (a_p)_n } g_n(z),\qquad n\in\mathbb{Z}_+. 
\end{equation}
Recall that a $R_I$-type continuos fraction is associated with a system of monic polynomials $\{ P_n(z) \}_{n=0}^\infty$, generated by
(\cite[p. 5]{cit_5100_Ismail_Masson_1995})
\begin{equation}
\label{f1_9} 
P_n(z) = (z - \mathbf{c}_n) P_{n-1}(z) - \lambda_n (z - \mathbf{a}_n) P_{n-2}(z),\qquad n=1,2,...,
\end{equation}
where $P_{-1}(z) := 0$, $P_0(z) := 1$, and
\begin{equation}
\label{f1_10} 
\lambda_{n+1}\not= 0,\quad P_n(\mathbf{a}_n)\not= 0.
\end{equation}
Polynomials $\{ P_n(z) \}_{n=0}^\infty$ are related to biorthogonal rational functions~\cite[Theorem 2.1]{cit_5100_Ismail_Masson_1995}.
The case $\mathbf{a}_n = 0$, $n\geq 2$, is related to general $T$-fractions~\cite{cit_1000_H_v_R__1986}.
It turnes out that this is the case for the monic polynomials $\{ G_n(z) \}_{n=0}^\infty$.
On the other hand, recall the following definition from~\cite{cit_95000_Z}:

\begin{definition}
\label{d1_1}
A set $\Theta = \left(
J_3, J_5, \alpha, \beta
\right)$,
where $\alpha>0$, $\beta\in\mathbb{R}$, $J_3$ is a Jacobi matrix and
$J_5$ is a semi-infinite real symmetric five-diagonal matrix with positive numbers on the second subdiagonal,
is said to be
\textbf{a Jacobi-type pencil (of matrices)}.
\end{definition}
Matrices $J_3$ and $J_5$ have the following form:
\begin{equation}
\label{ff1_5}
J_3 =
\left(
\begin{array}{cccccc}
b_0 & a_0 & 0 & 0 & 0 & \cdots\\
a_0 & b_1 & a_1 & 0 & 0 & \cdots\\
0 & a_1 & b_2 & a_2 & 0 &\cdots\\ 
\vdots & \vdots & \vdots & \ddots \end{array}
\right),\qquad a_k>0,\ b_k\in\mathbb{R},\ k\in\mathbb{Z}_+;
\end{equation}

\begin{equation}
\label{ff1_10}
J_5 =
\left(
\begin{array}{ccccccc}
\alpha_0 & \beta_0 & \gamma_0 & 0 & 0 & 0 & \cdots\\
\beta_0 & \alpha_1 & \beta_1 & \gamma_1 & 0 & 0 & \cdots\\
\gamma_0 & \beta_1 & \alpha_2 & \beta_2 & \gamma_2 & 0 & \cdots\\ 
0 & \gamma_1 & \beta_2 & \alpha_3 & \beta_3 & \gamma_3 & \cdots \\
\vdots & \vdots & \vdots &\vdots & \ddots \end{array}
\right),\ \alpha_n,\beta_n\in\mathbb{R},\ \gamma_n>0,\ n\in\mathbb{Z}_+.
\end{equation}

With a Jacobi-type pencil of matrices $\Theta$ one associates a system of polynomials
$\{ p_n(\lambda) \}_{n=0}^\infty$, such that
\begin{equation}
\label{ff1_15}
p_0(\lambda) = 1,\quad p_1(\lambda) = \alpha\lambda + \beta,
\end{equation}
and
\begin{equation}
\label{ff1_20}
(J_5 - \lambda J_3) \vec p(\lambda) = 0,
\end{equation}
where $\vec p(\lambda) = (p_0(\lambda), p_1(\lambda), p_2(\lambda),\cdots)^T$. Here the superscript $T$ means the transposition.
Polynomials $\{ p_n(\lambda) \}_{n=0}^\infty$ are said to be associated to the Jacobi-type pencil of matrices $\Theta$.
One can rewrite relation~(\ref{ff1_20}) in the scalar form:
$$ \gamma_{n-2} p_{n-2}(\lambda) + (\beta_{n-1}-\lambda a_{n-1}) p_{n-1}(\lambda) + (\alpha_n-\lambda b_n) p_n(\lambda) +
$$
\begin{equation}
\label{ff1_30}
+ (\beta_n-\lambda a_n) p_{n+1}(\lambda) + \gamma_n p_{n+2}(\lambda) = 0,\qquad n\in\mathbb{Z}_+,
\end{equation}
where $p_{-2}(\lambda) = p_{-1}(\lambda) = 0$, $\gamma_{-2} = \gamma_{-1} = a_{-1} = \beta_{-1} = 0$.

In the case of positive parameters $a_j,b_l$, the polynomials $\{ g_n(z) \}_{n=0}^\infty$ are connected with some
Jacobi type pencils and their associated polynomials $\{ p_n(\lambda) \}_{n=0}^\infty$. This connection resembles the connection 
between orthogonal polynomials
on the unit circle (OPUC) and orthogonal polynomials on $[-1,1]$.

In Section~2, the announced recurrence relation, a differential equation for $g_n(z)$, and Sobolev orthogonality
relations for $g_n(z)$ will be given in Theorem~1. Observe that the case $p=q=0$, leads to the exponential function.
The corresponding partial sums appeared in~\cite{cit_85000_Zagorodnyuk_CMA_2020}, as a particular case with $\rho=1$.

Of course, polynomials $g_n(z)$ have nice expressions for their coefficients. However, it is of interest to
get integral representations for $g_n(z)$, involving special functions. In particular, such integral representations would be useful
for obtaining various estimates, as well as in the Fourier series analysis.
We shall give two integral representations for $g_n(z)$ in Theorem~2.
Asymptotic properties of $g_n(z)$ and location of their zeros are also described by this theorem.
We shall discuss the partial sums of arbitrary power series with non-zero coefficients.
They are also related to biorthogonal rational functions.
Finally, we shall obtain a relation of polynomials $g_n(z)$ to Jacobi-type pencils and their associated polynomials.

\noindent
{\bf Notations. }
As usual, we denote by $\mathbb{R}, \mathbb{C}, \mathbb{N}, \mathbb{Z}, \mathbb{Z}_+$,
the sets of real numbers, complex numbers, positive integers, integers and non-negative integers,
respectively. 
For $k,l\in\mathbb{Z}$, we set $\mathbb{Z}_{k,l} := \{ j\in\mathbb{Z}: k\leq j\leq l\}$.
Set $\mathbb{T} := \{ z\in\mathbb{C}:\ |z|=1 \}$,
$\mathbb{D} := \{ z\in\mathbb{C}:\ |z|<1 \}$,
$\mathbb{D}_e := \{ z\in\mathbb{C}:\ |z|>1 \}$. 
By $\mathfrak{B}(\mathbb{T})$ we mean the set of all Borel subsets of $\mathbb{T}$.
By $\mathbb{P}$ we denote the set of all polynomials with complex coefficients.
For a complex number $c$ we denote
$(c)_0 = 1$, $(c)_k = c(c+1)...(c+k-1)$, $k\in\mathbb{N}$ (\textit{the shifted factorial or Pochhammer symbol}).
The generalized hypergeometric function is denoted by
$$ {}_p F_q(a_1,...,a_p; b_1,...,b_q;z) = \sum_{k=0}^\infty \frac{(a_1)_k ... (a_p)_k}{(b_1)_k ... (b_q)_k} \frac{z^k}{k!}, $$
where $p,q\in\mathbb{Z}_+$, $a_j,b_l\in\mathbb{C}$.

\section{The partial sums of the hypergeometric series and different kinds of orthogonality.}

Denote by $\mu_0$ the (probability) normalized arc length measure on $\mathbb{T}$,
which may be identified with the Lebesgue measure on $[0,2\pi)$.
We shall use the ideas from~\cite{cit_80000_Zagorodnyuk_JAT_2020} to obtain the following theorem.

\begin{theorem}
\label{t2_1}
Let $p,q\in\mathbb{Z}_+$, be some fixed numbers, and
$a_1,...,a_p$; $b_1,...,b_q$, be some parameters from $\mathbb{C}\backslash\{ 0,-1,-2,... \}$
(the case $p=0$ ($q=0$) means that $a_j$s (respectively $b_ls)$ are absent).
The following statements hold:

\begin{itemize}

\item[(a)]
Polynomials $g_n(z)$ from~(\ref{f1_5}) satisfy the following recurrence relation:
$$ \frac{ (n+1) (b_1+n) ... (b_q+n) }{ (a_1+n) ... (a_p+n) } ( g_{n+1}(z) - g_n(z) ) = $$
\begin{equation}
\label{f2_9} 
= z ( g_{n}(z) - g_{n-1}(z) ),\qquad n\in\mathbb{Z}_+,\quad g_{-1}(z) := 0. 
\end{equation}

\item[(b)]
The corresponding monic polynomials $G_n(z)$ from~(\ref{f1_7}) satisfy the following recurrence relation:
\begin{equation}
\label{f2_15} 
G_n(z) = (z+\delta_n) G_{n-1}(z) - \delta_{n-1} z G_{n-2}(z),\qquad n=1,2,..., 
\end{equation}
where $G_{-1}(z) := 0$, and $\delta_0 := 0$,
\begin{equation}
\label{f2_17} 
\delta_k := \frac{ k(b_1+k-1)...(b_q+k-1) }{ (a_1+k-1)...(a_p+k-1) },\qquad k\in\mathbb{N}.
\end{equation}
Therefore $G_n(z)$ are related to general $T$-fractions and biorthogonal rational functions.

\item[(c)]
Polynomials $g_n(z)$ obey the following differential equation:
\begin{equation}
\label{f2_19} 
\theta R g_n(z) - n R g_n(z) = 0,\qquad n\in\mathbb{Z}_+,
\end{equation}
where 
\begin{equation}
\label{f2_27} 
R  = \frac{d}{dz} \prod_{j=1}^q (\theta + b_j -1) - \prod_{j=1}^p (a_j + \theta),\qquad \theta := z\frac{d}{dz}.
\end{equation}

\item[(d)] Polynomials $g_n(z)$ are Sobolev orthogonal polynomials on the unit circle:
\begin{equation}
\label{f2_30}
\int_{\mathbb{T}} \left( g_n(z), g_n'(z),..., g_n^{(\rho)}(z) \right) M \overline{
\left( \begin{array}{cccc} g_m(z) \\
g_m'(z) \\
\vdots \\
g_m^{(\rho)}(z) \end{array} \right)
} 
d\mu_0 = \delta_{n,m},\qquad n,m\in\mathbb{Z}_+,
\end{equation}
where
\begin{equation}
\label{f2_35}
M = (c_0(z),c_1(z),...,c_\rho(z))^T (\overline{c_0(z)}, \overline{c_1(z)},...,\overline{c_\rho(z)}),
\end{equation}
and $c_j(z)\in\mathbb{P}$, are the coefficients of the differential operator:
\begin{equation}
\label{f2_38}
-\frac{ n! (b_1)_n ... (b_q)_n }{ (a_1)_{n+1} ... (a_p)_{n+1} } R = \sum_{l=0}^\rho c_l(z) \frac{ d^l }{ dz^l},
\end{equation}
with $\rho = \max(p,q+1)$.

\end{itemize}

Note that the case $p=0$ ($q=0$) means that all $(a_j)_k$, $(a_j+k)$ (respectively $(b_l)_k$, $(b_l+k)$) in
the above formulas are replaced by $1$. The same takes place with the products $\prod_{j=1}^p$ and $\prod_{j=1}^q$.
\end{theorem}
\textbf{Proof.}

\noindent
$(a)$: Observe that
$$ g_{n}(z) - g_{n-1}(z) = \frac{ (a_1)_n ... (a_p)_n }{ (b_1)_n ... (b_q)_n } \frac{z^n}{n!},\qquad n\in\mathbb{Z}_+, $$
where $g_{-1} = 0$.
Using this relation with
\begin{equation}
\label{f2_40}
z^{n+1} = z z^n,\qquad n\in\mathbb{Z}_+, 
\end{equation}
we immediately obtain the required recurrence relation~(\ref{f2_9}).

\noindent
$(b)$: It follows directly from~$(a)$.

\noindent
$(c)$: Using the known idea of proof for the differential equation of ${}_p F_q$ (\cite{cit_5150_Rainville}) we may write:
$$ \theta\prod_{j=1}^q (\theta+b_j-1) g_n(z) = 
\sum_{k=0}^n \frac{ (a_1)_k ... (a_p)_k }{ (b_1)_k ... (b_q)_k } \frac{1}{k!} 
k \prod_{j=1}^q (k+b_j-1) z^k = $$
$$ = z \sum_{k=1}^n \frac{ (a_1)_k ... (a_p)_k }{ (b_1)_{k-1} ... (b_q)_{k-1} } \frac{z^{k-1}}{(k-1)!} = 
z \sum_{l=0}^{n-1} \frac{ (a_1)_{l+1} ... (a_p)_{l+1} }{ (b_1)_{l} ... (b_q)_{l} } \frac{ z^l }{ l! } = $$
$$ = z \sum_{l=0}^{n-1} \frac{ (a_1)_{l} ... (a_p)_{l} }{ (b_1)_{l} ... (b_q)_{l} } \frac{1}{l!} \prod_{j=1}^p (a_j+l)
z^l = z \prod_{j=1}^p (a_j+\theta) \sum_{l=0}^{n-1} \frac{ (a_1)_{l} ... (a_p)_{l} }{ (b_1)_{l} ... (b_q)_{l} } \frac{z^l}{l!} = $$
$$ = z \prod_{j=1}^p (a_j+\theta) \left(
g_n(z) - \frac{ (a_1)_n ... (a_p)_n }{ (b_1)_n ... (b_q)_n } \frac{z^n}{n!}
\right),\qquad n\in\mathbb{Z}_+. $$
Then
\begin{equation}
\label{f2_42}
-\frac{ n! (b_1)_n ... (b_q)_n }{ (a_1)_{n+1} ... (a_p)_{n+1} } R g_n(z) = z^n,\qquad n\in\mathbb{Z}_+,
\end{equation}
where $R$ is defined by~(\ref{f2_27}).
Here we assumed that $p,q\in\mathbb{N}$, while the other cases are similar and lead to the same formula. 
Using 
$$ \theta z^n = n z^n,\qquad n\in\mathbb{Z}_+, $$
and relation~(\ref{f2_42}) we obtain the differential equation~(\ref{f2_19}).

\noindent
$(d)$: It follows from the orthonormality relations for $\{ z^n \}_{n=0}^\infty$, and relations~(\ref{f2_42}),(\ref{f2_38}).
$\Box$

Integral representations for polynomials $g_n(z)$ and some other their basic properties are described in
the next theorem.

\begin{theorem}
\label{t2_2}
In conditions of Theorem~\ref{t2_1}
the following statements hold:

\begin{itemize}

\item[(i)]
If $p\leq q$, then polynomials $g_n$ admit the following integral representation:
$$ g_n(e^{i\tau}) = 
\frac{1}{2\pi} \int_0^{2\pi} \left( \frac{ 1 - e^{ i(n+1)(\tau - t) } }{ 1 - e^{ i(\tau - t) } } \right)
{}_p F_q(a_1,...,a_p; b_1,...,b_q;e^{it}) dt, $$
\begin{equation}
\label{f2_45} 
\tau\in [0,2\pi),\ n\in\mathbb{Z}_+.
\end{equation}

\item[(ii)]
Polynomials $g_n$ have the following integral representation:
$$ g_n(x) = $$
$$ = -(n+1) x^{n+1} \int_{-\infty}^x t^{-n-2} {}_{p+1} F_{q+1} (-n,a_1,...,a_p; -n-1, b_1,...,b_q;t) dt, $$
\begin{equation}
\label{f2_47} 
x<0,\ n\in\mathbb{Z}_+.
\end{equation}

\item[(iii)]
Polynomials $g_n$ have simple roots.
If $p\leq q$, and
\begin{equation}
\label{f2_48} 
0 < a_j \leq b_j,\qquad j\in\mathbb{Z}_{1,p};
\end{equation}
\begin{equation}
\label{f2_49} 
b_k\geq 1,\qquad k\in\mathbb{Z}_{p+1,q},
\end{equation}
then all roots of $g_n$ are located in $\mathbb{D}_e\backslash(1,\infty)$.

\item[(iv)] Polynomials $g_n(z)$ tend to ${}_p F_q(a_1,...,a_p; b_1,...,b_q;z)$, as $n\rightarrow\infty$, in $\mathbb{T}$
(in $\mathbb{C}$), if $p=q+1$ (respectively $p\leq q$).

\end{itemize}

\end{theorem}
\textbf{Proof.}

\noindent
$(i)$: If $p\leq q$, then the function $g(z):= {}_p F_q(a_1,...,a_p; b_1,...,b_q;z)$, is analytic in the whole plane.
In particular, we have
\begin{equation}
\label{f2_52} 
g(e^{i\tau}) = \sum_{k=0}^\infty \xi_k e^{ik\tau} = \lim\limits_{n\to\infty} g_n(e^{i\tau}),\qquad \tau\in[0,2\pi), 
\end{equation}
where 
\begin{equation}
\label{f2_53} 
\xi_k := \frac{ (a_1)_k ... (a_p)_k }{ (b_1)_k ... (b_q)_k } \frac{ 1 }{ k! }.
\end{equation}
Let us check that
\begin{equation}
\label{f2_54} 
\frac{1}{2\pi} \int_0^{2\pi} g(e^{i\tau}) e^{-ij\tau} d\tau = \xi_j,\qquad j\in\mathbb{Z}_+.
\end{equation}
Observe that
\begin{equation}
\label{f2_56} 
\left|
g_n(e^{i\tau}) e^{-ij\tau}
\right|
\leq \sum_{k=0}^n |\xi_k| \leq \sum_{k=0}^\infty |\xi_k| =: C < \infty. 
\end{equation}
By~(\ref{f2_52}), (\ref{f2_56}) and the dominated convergence theorem it follows that
\begin{equation}
\label{f2_58} 
\frac{1}{2\pi} \int_0^{2\pi} g_n (e^{i\tau}) e^{-ij\tau} d\tau \rightarrow_{n\rightarrow\infty}
\frac{1}{2\pi} \int_0^{2\pi} g(e^{i\tau}) e^{-ij\tau} d\tau.
\end{equation}
The left-hand side of~(\ref{f2_58}) is equal to $\xi_j$, when $n\geq j$. Therefore relation~(\ref{f2_54}) hols true.
We may write
\begin{equation}
\label{f2_60} 
g_n(e^{i\tau}) = \sum_{k=0}^n \xi_k e^{ik\tau} = 
\frac{1}{2\pi} \int_0^{2\pi} g(e^{it}) 
\left(
\sum_{k=0}^n e^{ik (\tau - t)}
\right) dt,\quad \tau\in[0,2\pi). 
\end{equation}
Notice that the following relation:
\begin{equation}
\label{f2_62} 
d_n = d_n(u) := \sum_{k=0}^n u^k = \frac{1-u^{n+1}}{1-u},\qquad n\in\mathbb{Z}_+, 
\end{equation}
holds for all complex $u$ (not only for $u\in\mathbb{D}$). It follows from the recurrence relation:
$d_{n+1} - u d_n = 1$, with $d_0 = 1$.
By~(\ref{f2_60}) and (\ref{f2_62}) we obtain relation~(\ref{f2_45}).

\noindent
$(ii)$: Formula~(\ref{f2_47}) can be checked by the direct integration of the polynomial under the integral sign,
and some algebraic simplifications.

\noindent
$(iii)$: 
The simplicity of zeros of $g_n$ follows from the three-term recurrence relation~(\ref{f2_9}) (notice that $g_n(0)=1$).

Let $p\leq q$, and conditions~(\ref{f2_48}),(\ref{f2_48}) be satisfied.
Since coefficients of $g_n$ are positive, it takes positive values on $(0,+\infty)$.
Conditions~(\ref{f2_48}),(\ref{f2_49}) imply that
the coefficients of $g_n$ form a monotone sequence. Thus, one can apply the Enestr\"om--Kakeya Theorem 
(\cite[p. 136]{cit_5250_Marden}) to conclude that the roots of $g_n$ are located in $\mathbb{D}_e$.

\noindent
$(iv)$: It is a known property of ${}_p F_q$.
$\Box$

\noindent
\textbf{Generalizations.}
Consider a power series
\begin{equation}
\label{f2_64} 
f(z) = \sum_{k=0}^\infty d_k z^k,\qquad d_k\in\mathbb{C}\backslash\{ 0 \}.
\end{equation}
Define
\begin{equation}
\label{f2_67} 
f_n (z) = \sum_{k=0}^n d_k z^k,\qquad n\in\mathbb{Z}_+,
\end{equation}
\begin{equation}
\label{f2_70} 
F_n (z) = \frac{1}{d_n} f_n(z),\qquad n\in\mathbb{Z}_+.
\end{equation}
We have 
\begin{equation}
\label{f2_72} 
\frac{1}{d_n} (f_{n}(z) - f_{n-1}(z)) = z^n,,\qquad n\in\mathbb{Z}_+, 
\end{equation}
where $f_{-1}:=0$.
By~(\ref{f2_40}),(\ref{f2_72}) we obtain a recurrence relation for $f_n$. For the monic polynomials
$F_n$ it takes the following form:
\begin{equation}
\label{f2_74} 
F_n(z) = \left(
z + \frac{d_{n-1}}{d_n}
\right) 
F_{n-1}(z) - 
\frac{ d_{n-2} }{ d_{n-1} }
z F_{n-2}(z),\qquad n=1,2,..., 
\end{equation}
where $F_{-1}(z) := 0$, and $d_{-1} := 1$.
Consequently, polynomials $F_n(z)$ are also related to biorthogonal rational functions.
In order to obtain some equations for $f_n(z)$ with respect to $z$ (e.g., differential, difference, $q$-difference equations),
it looks promising to consider the power series $f(z)$, corresponding to the basic hypergeometric function,
elliptic hypergeometric functions and other special functions.

Suppose now that all coefficients $d_k$ are positive.
We may write:
\begin{equation}
\label{f2_76} 
\mathop{\rm Re}\nolimits f_n(e^{i\tau}) = \sum_{k=0}^n d_k T_k(x), 
\end{equation}
$$ \mathop{\rm Im}\nolimits  f_{n+1}(e^{i\tau}) = \sin \tau \sum_{k=1}^{n+1} d_k U_{k-1}(x) =
\sin \tau \sum_{j=0}^n d_{j+1} U_j(x), $$
\begin{equation}
\label{f2_78} 
x=\cos\tau,\ \tau\in(0,\pi),\ n\in\mathbb{Z}_+, 
\end{equation}
where $T_k(x) = \cos(k\arccos x)$,
$U_k(x) = \frac{ \sin((k+1)\arccos x) }{ \sqrt{1-x^2} }$, are Chebyshev polynomials of the first and the second kinds.
Therefore the following systems of polynomials:
$$\left\{ 
\mathop{\rm Re}\nolimits f_n(e^{i\arccos x}) 
\right\}_{n=0}^\infty,\quad
\left\{ \frac{1}{\sin \arccos x} 
\mathop{\rm Im}\nolimits f_{n+1}(e^{i\arccos x}) 
\right\}_{n=0}^\infty, $$
are \textit{modified kernel polynomials}, see~\cite[formula (5)]{cit_87000_Zagorodnyuk_Modified_kernel_polynomials}.
They are associated to Jacobi-type pencils and posess special orthogonality relations.

\begin{center}
{\large\bf 
On the generalized hypergeometric function, Sobolev orthogonal polynomials and biorthogonal rational functions.}
\end{center}
\begin{center}
{\bf S.M. Zagorodnyuk}
\end{center}

It turned out that the partial sums
$g_n(z) = \sum_{k=0}^n \frac{(a_1)_k ... (a_p)_k}{(b_1)_k ... (b_q)_k} \frac{z^k}{k!}$,
of the generalized hypergeometric series ${}_p F_q(a_1,...,a_p; b_1,...,b_q;z)$, with
parameters $a_j,b_l\in\mathbb{C}\backslash\{ 0,-1,-2,... \}$, are Sobolev orthogonal polynomials.
The corresponding monic polynomials $G_n(z)$ are polynomials of $R_I$ type, and therefore
they are related to biorthogonal rational functions.
Polynomials $g_n$ possess a differential equation (in $z$), and a recurrence relation (in $n$).
We study integral representations for $g_n$, and some other their basic properties.
Partial sums of arbitrary power series with non-zero coefficients are shown to be also related to biorthogonal rational functions.
We obtain a relation of polynomials $g_n(z)$ to Jacobi-type pencils and their associated polynomials.

\vspace{1.5cm} 

V. N. Karazin Kharkiv National University \newline\indent
School of Mathematics and Computer Sciences \newline\indent
Department of Higher Mathematics and Informatics \newline\indent
Svobody Square 4, 61022, Kharkiv, Ukraine

Sergey.M.Zagorodnyuk@gmail.com; Sergey.M.Zagorodnyuk@univer.kharkov.ua

}
\end{document}